

\input amstex

\documentstyle{amsppt}

\loadbold

\magnification=\magstep1

\pageheight{9.0truein}
\pagewidth{6.5truein}


\def\k{\overline{k}}
\def\nxn{n{\times}n}

\def\ker{\operatorname{ker}}

\def\Rel{\operatorname{Rel}}
\def\ann{\operatorname{ann}}
\def\nxn{n{\times}n}
\def\sgn{\operatorname{sgn}}
\def\trace{\operatorname{trace}}
\def\Trace{\operatorname{Trace}}
\def\x{\bold{x}}
\def\y{\bold{y}}
\def\simplen{\operatorname{simple}_n}
\def\ZZ{\Bbb{Z}}

\def\Art{1}
\def\BecWei{2}
\def\Bok{3}
\def\For{4}
\def\Let{5}
\def\McCRob{6}
\def\Pap{7}
\def\Pro{8}
\def\Raz{9}

\topmatter

\title Counting equivalence classes of irreducible representations \endtitle

\date June 2001. \enddate

\rightheadtext{Counting irreducible representations}

\author Edward S. Letzter \endauthor

\abstract Let $n$ be a positive integer, and let $R$ be a (possibly
infinite dimensional) finitely presented algebra over a computable
field of characteristic zero. We describe an algorithm for deciding
(in principle) whether $R$ has at most finitely many equivalence
classes of $n$-dimensional irreducible representations. When $R$ does
have only finitely many such equivalence classes, they can be
effectively counted (assuming that $k[x]$ posesses a factoring
algorithm).  \endabstract

\address Department of Mathematics, Temple University, Philadelphia,
PA 19122 \endaddress

\email letzter\@math.temple.edu\endemail

\thanks The author's research was supported in part by NSF grant DMS-9970413.
Also, a part of this research was completed while the author was a participant
(February 2000) in the MSRI program on noncommutative algebra. \endthanks

\endtopmatter

\document

\head 1. Introduction \endhead

Let $n$ be a positive integer, fixed throughout. In \cite{\Let} we
observed that the existence of $n$-dimensional irreducible
representations of finitely presented noncommutative algebras can be
algorithmically decided. In this note we outline a procedure for
effectively ``counting'' the number of such irreducible
representations, up to equivalence, in characteristic zero. Our
approach combines standard computational commutative algebra with
results from \cite{\Art} and \cite{\Raz}.

\subhead 1.1 \endsubhead Assume that $k$ is a computable field of
characteristic zero, and that $\k$ is the algebraic closure of
$k$. 

Henceforth, let
$$R = k\{ X_1,\ldots, X_s \}/\langle f_1,\ldots, f_t \rangle,$$
for some fixed choice of $f_1,\ldots , f_t$ in the free associative
$k$-algebra $k\{X_1,\ldots , X_s \}$. In a slight abuse of notation,
``$X_\ell$'' will also denote its image in $R$, for $1 \leq \ell \leq
s$.

By an {\sl $n$-dimensional representation of $R$\/} we will always
mean a unital $k$-algebra homomorphism from $R$ into the $k$-algebra
$M_n(\k)$ of $\nxn$ matrices over $\k$. Representations $\rho, \rho'
\colon R \rightarrow M_n(\k)$ are {\sl equivalent\/} if there exists a
matrix $Q \in GL_n(\k)$ such that
$$\rho'(X) = Q\rho(X) Q^{-1},$$
for all $X \in R$. 

We will say that the representation $\rho\colon R \rightarrow M_n(\k)$
is {\sl irreducible\/} when $\k \rho (R) = M_n(\k)$ (cf\. \cite{1, \S
9}). Observe that $\rho$ is irreducible if and only if $\rho\otimes 1
\colon R \otimes_k \k \rightarrow M_n(\k)$ is surjective, if and only
if $\rho \otimes 1$ is irreducible in the more common use of the term.
(In particular, our approach below will use calculations over the
computable field $k$ to study representations over the algebraically
closed field $\k$.)

\subhead 1.2 \endsubhead The existence of an $n$-dimensional
representation of $R$ depends only on the consistency of a system of
algebraic equations, over $k$, in $(t.n^2)$-many
variables. Consequently, the existence of $n$-dimensional
representations of $R$ is decidable (in principle) using Groebner
basis methods. This idea is extended in \cite{\Let} to give a
procedure for deciding the existence of $n$-dimensional irreducible
representations. On the other hand, posessing a nonzero finite
dimensional representation is a Markov property, and so the existence
-- in general -- of a finite dimensional representation of $R$ cannot
be effectively decided, by \cite{\Bok}.

We now state our main result; the proof will be presented in \S 2.

\proclaim{Theorem} Having at most most finitely many equivalence classes of
irreducible $n$-\-dimen\-sional representations is an algorithmically
decidable property of $R$. \endproclaim

\subhead 1.3 \endsubhead Assume that $k[x]$ is equipped with a factoring
algorithm. If it has been determined that $R$ has at most finitely many
equivalence classes of $n$-dimensional irreducible representations, these
equivalence classes can (in principle) be effectively counted; see (2.9).

\head 2. Proof of Theorem \endhead

\subhead 2.1 \endsubhead (i) Set
$$B = k [ x_{ij}(\ell) : 1 \leq i,j \leq n, 1 \leq \ell \leq s ].$$
For $1 \leq \ell \leq s$, let $\x _\ell$ denote the $\nxn$ generic
matrix $(x_{ij}(\ell))$, in $M_n(B)$. For $g \in k\{ X_1, \ldots , X_s
\}$, let $g(\x)$ denote the image of $g$, in $M_n(B)$, under the
canonical map
$$k\{ X_1,\ldots , X_s \} @> \; X_\ell \longmapsto \x_\ell \; >>
M_n(B) .$$
Identify $B$ with the center of $M_n(B)$.

(ii) Let $\Rel (M_n(B))$ be the ideal of $M_n(B)$ generated by
$f_1(\x), \ldots , f_t (\x)$.

(iii) Let $\Rel (B)$ denote the ideal of $B$ generated by the entries
of the matrices $f_1(\x)$, $\ldots$, $f_t(\x)$ $\in$ $M_n(B)$. Note that
$$\Rel(B) = \Rel(M_n(B)) \cap B .$$

(iv) Let
$$A = k\{ \x_1, \ldots , \x_s \},$$
the $k$-subalgebra of $M_n(B)$ generated by the generic matrices
$\x_1, \ldots , \x_s$.  Set
$$\Rel (A) \; = \; \Rel (M_n(B)) \cap A .$$

\subhead 2.2 \endsubhead Every $n$-dimensional representation of $R$ can be
written in the form
$$R @> X_\ell \longmapsto \x_\ell + \Rel(A) >> \left(\frac{A}{\Rel(A)}
\right) @> \text{inclusion} >> \left(\frac{M_n(B)}{\Rel
(M_n(B))}\right) \longrightarrow M_n(\k) , $$
and every $k$-algebra homomorphism
$$M_n(B) / \Rel(M_n(B)) \; \rightarrow \; M_n(\k)$$
is completely determined by the induced map
$$B/\Rel(B) \rightarrow \k .$$
For each representation $\rho \colon R \rightarrow M_n(\k)$, let
$\chi_\rho \colon B \rightarrow \k$ be the homomorphism (with $\Rel(B)
\subseteq \ker\chi_\rho$) given by this correspondence.

\subhead 2.3 \endsubhead Let $T$ be the $k$-subalgebra of $B$
generated by the coefficients of the characteristic polynomials of
elements in $A$. (Since the characteristic of $k$ is zero, $T$ is in
fact generated by the traces, as $\nxn$ matrices, of the elements in
$A$.)  Set
$$\Rel(T) = \Rel (B) \cap T .$$
Note, when $\rho,\rho' \colon R \rightarrow M_n(\k)$ are equivalent
representations, that the restrictions of $\chi_\rho$ and
$\chi_{\rho'}$ to $T$ will coincide.

\subhead 2.4 \endsubhead Let $\simplen (R)$ denote the set of
equivalence classes of irreducible $n$-dimensional representations of
$R$. By (2.3) there is a well-defined function
$$\Phi : \simplen (R) \; \longrightarrow \; V(\Rel (T)),$$
where $V(\Rel(T))$ denotes the $\k$-affine algebraic set of points on
which the polynomials in $\Rel(T)$ vanish. It follows from \cite{\Art,
pp\. 558--559} that $\Phi$ is injective.

\subhead 2.5 \endsubhead (i) Recall the {\sl $m$th standard
identity\/}
$$s_m \; = \; \sum_{\sigma \in S_m} (\sgn \sigma)Y_{\sigma(1)}\cdots
Y_{\sigma(m)} \; \in \; \ZZ \{ Y_1,\ldots,Y_m \} .$$
If $\Lambda$ is a commutative ring, then the Amitsur-Levitzky Theorem ensures
that $M_n(\Lambda)$ satisfies $s_m$ if and only if $m \geq 2n$; see, for
example, \cite{\McCRob, 13.3.2, 13.3.3}.

(ii) Let $S$ denote the finite subset of $T$ ($\subseteq B$) comprised
of
$$\trace\left(M_0 \cdot s_{2(n-1)}(M_1,\ldots,M_{2(n-1)})\right) ,$$
for all monic monomials $M_0,\ldots,M_{2(n-1)}$, in the generic
matrices $\x_1,\ldots,\x_s$, of length less than
$$p = n\sqrt{2n^2/(n-1)+ 1/4} + n/2 -2 .$$
(The choice of $p$ will follow from \cite{\Pap}; see \cite{\Let,
2.2}.)  Let $\rho\colon R \rightarrow M_n(\k)$ be a representation. It
now follows from \cite{\Let, \S 2} that $\rho$ is irreducible if and
only if
$$S \; \not\subseteq \; \ker\chi_\rho .$$
(Other sets of polynomials can be substituted for $S$; see \cite{\Let,
2.6vi,vii}.)

\subhead 2.6 \endsubhead (i) Set
$$W = \; V(\Rel (T)) \setminus V(S) .$$
Combining (2.4) with (2.5ii), we obtain a bijection
$$\Phi : \simplen (R) \; \longrightarrow \; W .$$

(ii) Set
$$J = \ann_B\left(\frac{\Rel(B) + B.S}{\Rel(B)}\right), \quad
\text{and} \quad I = J\cap T = \ann_T\left(\frac{\Rel(T) +
T.S}{\Rel(T)}\right) .$$
A finite generating set for $J$ can be specified, using standard
methods, and we can identify $T/I$ with its image in $B/J$. Since
$V(I)$ is the Zariski closure of $W$, to prove the theorem it suffices
to find an effective procedure for determining whether or not $T/I$ is
finite dimensional. (When not indicated otherwise, ``dimension''
refers to ``dimension as a $k$-vector space.'')

\subhead 2.7 \endsubhead (i) For the generic matrices
$\x_1,\ldots,\x_s$, set $\Trace =$
$$\left\{ \trace ( \y_1 \y_2 \cdots \y_u ) : \text{$\y_1,\ldots,\y_u
\in \{ \x_1, \ldots , \x_s \}$ and $1 \leq u \leq n^2$} \right\} .$$
In \cite{\Raz} (cf\. \cite{\For, p\. 54}) it is shown that $T = k[
\Trace ]$. (A larger finite generating set for $T$ was established
in \cite{\Pro}.)

(ii) By (2.6ii), to prove the theorem it remains to find an algorithm
for deciding whether the monomials in $\Trace$ ($\subseteq B$) are
algebraic over $k$, modulo $J$. We accomplish this task using a
variant of the subring membership test (cf., e.g., \cite{\BecWei,
p\. 270}): Let $C$ be a commutative polynomial ring, over $k$, in $m$
variables. Let $L$ be an ideal -- equipped with an explicitly given
list of generators -- in $C$. Choose $f \in C$. Observe that $f$ is
algebraic over $k$, modulo $L$, if and only if $L \cap k[f] \ne
\{0\}$. Next, embed $C$, in the obvious way, as a subalgebra of the
polynomial ring $C' = k(t) \otimes _k C$. Observe that $L \cap k[f]
\ne \{0\}$ if and only if $1$ is contained in the ideal $(t-f).C' +
L.C'$ of $C'$. Hence, the decidability of ideal membership in $C'$
implies the decidability of algebraicity modulo $L$ in $C$.

The proof of the theorem follows.

\subhead 2.8 \endsubhead Roughly speaking, the complexity of the
procedure described in (2.1 -- 2.7) varies according to the degrees of
the polynomials involved in deciding the algebraicity of $\Trace$
modulo $J$. Note, for example, that the degrees of the members of $S$
can be as large as $p^{2n-1}$, for $p$ as in (2.5ii).

\subhead 2.9 \endsubhead Assume that it has already been determined
that the number (equal to $|W|$) of equivalence classes of irreducible
$n$-dimensional representations of $R$ is finite. Further assume that
$k[x]$ is equipped with a factoring algorithm. We conclude our study
by sketching a procedure for calculating -- in principal -- this
number.

Set $D = T/I$, and identify $D$ with the (finite dimensional)
$k$-subalgebra of $B/J$ generated by the image of $\Trace$. Since
$B/J$ can be given a specific finite presentation, finding a $k$-basis
$E$ for $D$ amounts to solving systems of polynomial equations in $B$,
and this task can be accomplished employing elimination methods. Next,
using the regular representation of $D$, and the finite presentation
of $B/J$, we can algorithmically specify $E$ as a set of commuting
$m{\times}m$ matrices over $k$, for some $m$. Furthermore, the
nilradical $N(D)$ will be precisely the set of elements of $D$ whose
traces, as $m{\times}m$ matrices, are zero. Consequently, we can
effectively compute the dimension of $D/N(D)$. This dimension is equal
to $|W|$.

\enddocument